%
%
%
%
%
\input amstex
\documentstyle{amsppt}
\def\curraddr#1\endcurraddr{\address
        {\it Current address \/}: #1\endaddress}
 
\topmatter
\title On the removable singularities\\
for meromorphic mappings\endtitle
\rightheadtext{ON THE REMOVABLE SINGULARITIES}
\author E.M.Chirka\endauthor
\address Steklov Mathematical Institute, Vavilov st. 42, Moscow
GSP-1, 117966 Russia\endaddress
\email chirka\@mph.mian.su\endemail
 
 
\abstract
If $E$ is a nonempty closed subset of the locally finite Hausdorff
(2n-2)-measure on an n-dimensional complex manifold $\Omega$ and all points
of $E$ are nonremovable for a meromorphic mapping of $\Omega \setminus E$
into a compact K\"ahler manifold, then $E$ is a pure (n-1)-dimensional
complex analytic subset of $\Omega$.
\endabstract
 
\thanks The final version of this paper will be
submitted for publication elsewhere
\endthanks
\endtopmatter
 
\document
\subhead 1\endsubhead
This paper was inspired by the following question of E.L.Stout \cite{7}:
Let $E$ be a closed subset of the complex projective space
${\Bbb P}^n \; (n \ge 2)$ such that the Hausdorff $(2n-2)$-measure of $E$
(with respect to the Fubini-Study metric) is less then that of any complex
hyperplane of ${\Bbb P}^n$. If it is true that $E$ is then a removable
singularity for meromorphic functions? Using natural projections of
${\Bbb P}^n$ onto hyperplanes G.Lupacciolu \cite{5} has shown the removability
of $E$ under additional conditions on the sizes of $E$ and a maximal ball
in the complement. (The projection of ${\Bbb P}^n$
onto a hyperplane does not
decrease Hausdorff measures, as it take place in the Euclidean space.)
Using an Oka--Nishino theorem \cite{6} on pseudoconcave sets we prove
here the following.
 
\proclaim{Theorem} Let $E$  be a closed subset of the locally finite
Hausdorff $(2n-2)$-measure on an $n$-dimensional complex manifold
$\Omega $  and let $f$  be a meromorphic mapping of
$\Omega \setminus E$  into a complex manifold $X$.
If $X$  has the meromorphic extension property and $E$
does not contain any $(n-1)$-dimensional closed analytic subset of
$\Omega $ then $f$ is continued to a meromorphic mapping of
$\Omega $ into $X$.
\endproclaim
 
Here we say that $X$ has the meromorphic extension property, if any
meromorphic map $\varphi : T \to X$ of the "squeezed polydisc"
$$
T = {(z,w) \in {\Bbb C}^{n-1}_{z} \times {\Bbb C}_{w} : |z| < r,
|w| < 1  \text{or} |z| < 1, 1-r < |w| < 1 },
$$
$0 < r < 1, \ n \ge 2,$ extends to a meromorphic map
$\tilde {\varphi} : U \to X$ of the unit polydisc
$U: |z| < 1, |w| < 1.$  By a recent result of
S.M.Ivashkovich \cite{3} every compact K\"ahler manifold $X$ has the
meromorphic extension property, so we have a lot of nice examples of such
$X$. The case of meromorphic functions ($X = {\Bbb P}^1$) is almost trivial
in the consideration: every meromorphic function in a squeezed polydisc is
represented as a ratio of two holomorphic functions (see \cite{4}) and
thus it is meromorphically continued into $U$. So the answer on the question
of Stout is positive because the mentioned $E$ can not contain any complex
analytic subset of ${\Bbb P}^n$ (by a Chow theorem such a set is algebraic,
and the Hausdorff $(2n-2)$-measure of it is not less then the measure of
a complex hyperplane (see, e.g. \cite{C}). As pointed me E.L.Stout, this
question was solved already by his student Mark Lawrence by a similar method.
 
\subhead 2\endsubhead
A closed subset $\Sigma $ of a complex manifold is called locally
pseudoconcave if for every point $a \in {\Sigma }$ there exists a Stein
neighbourhood $V$ such that $V \setminus \Sigma $ is Stein.
 
Let $E^{\prime}$ be the set of points $a \in E$ such that $f$ meromorphically
extends into a neighbourhood of $a$. Then $S:= E \setminus E^{\prime }$ is
closed. As the complement to $E$ is locally connected in $\Omega $, these
local meromorphic continuations of $f$ in points of $E^{\prime }$ glue
together into the unique meromorphic map
$f: \Omega \setminus S \to X$ (we preserve the notation $f$).
 
The proof of the Theorem is accomplished now in two steps. Firstly we
prove that $S$ is locally pseudoconcave (Lemma 1), and secondly we prove
that $S$ is complex analytic (Lemma 2).
 
\proclaim{Lemma 1} $S$ is locally pseudoconvex in $\Omega$.
\endproclaim
 
\demo{Proof} Let $a \in S,\ V \ni a$ is biholomorphic to a ball in
${\Bbb C}^n$ and $\varphi : T \to V \setminus S$ is a
holomorphic embedding. Then ($V$ is Stein) $\varphi $ extends to a
holomorphic embedding $\tilde {\varphi } : U \to V$ (see \cite{2}).
As $X$ has the meromorphic extension property, the meromorphic map
$f\circ \varphi : T \to X$ extends to a meromorphic map of
$U$ into $X$ and thus $f$ meromorphically continues into the domain
$\tilde {\varphi }(U) \subset V$. By the definition of
$S,\ \tilde {\varphi }(U)$ does not intersect $S$, and thus we have proved
that $V \setminus S$ sutisfies the condition "$p_{7}$-convexity" of
Docquier--Grauert \cite{2}. It follows that $V \setminus S$ is Stein.
\enddemo
 
\proclaim{Lemma 2}  Let $S$  be a nonempty locally pseudoconcave
subset of the finite Hausdorff $(2n-2)$-measure on an $n$-dimensional
complex manifold $\Omega$. Then $S$ is a pure $(n-1)$-dimensional complex
analytic subset of $\Omega$.
\endproclaim
 
\demo{Proof}  The statement is local, so we can assume that $\Omega $
is a domain in ${\Bbb C}^n$. Let $a \in S$ and coordinates
$(z,w),\ z = (z_{1},...,z_{n-1})$ in ${\Bbb C}^n$ are choosen in such a way
that $a = 0$ and the set $S \cap \{z = 0\}$ is finite or countable (it
can be done obviously). Then there exists a neighbourhood
$V = V^{\prime } \times V_{n} \ni 0$ such that the projection of $S \cap V$
into $V^{\prime }$ is proper. It follows that fibres
$S \cap V \cap \{z = c \}$ are finite for almost every $c \in V^{\prime }$.
By an Oka--Nishino theorem \cite{6} \ $S \cap V$ is then a complex analytic
set of pure dimension $n-1$.
\enddemo
 
I would like to thank G.Lupacciolu for the reprint of his paper \cite{5},
which stimulated this work.

\enddocument